\documentclass[]{interact}

\usepackage{epstopdf}
\usepackage[caption=false]{subfig}

\usepackage[numbers,sort&compress]{natbib}
\bibpunct[, ]{[}{]}{,}{n}{,}{,}
\makeatletter
\def\NAT@def@citea{\def\@citea{\NAT@separator}}
\makeatother
\usepackage{anyfontsize}
\theoremstyle{plain}
\newtheorem{theorem}{Theorem}[section]
\newtheorem{lemma}[theorem]{Lemma}
\newtheorem{corollary}[theorem]{Corollary}
\newtheorem{proposition}[theorem]{Proposition}

\theoremstyle{definition}
\newtheorem{definition}[theorem]{Definition}
\newtheorem{remark}[theorem]{Remark}

\usepackage{xcolor}

\theoremstyle{remark}

\begin{document}
	

	\title{$A$-approximate point spectrum of $A$-bounded operators in semi-Hilbertian spaces}
  \author{\name {Arup Majumdar \thanks{Arup Majumdar (corresponding author). Email address: arupmajumdar93@gmail.com}and P. Sam Johnson \thanks{ P. Sam Johnson. Email address: sam@nitk.edu.in}} \affil{Department of Mathematical and Computational Sciences, \\
			National Institute of Technology Karnataka, Surathkal, Mangaluru 575025, India.}}
  
	\maketitle

	\begin{abstract}
      This paper delves into several characterizations of $A$-approximate point spectrum of A-bounded operators acting on a complex semi-Hilbertian space $H$ and also investigates properties of the $A$-approximate point spectrum for the tensor product of two $A^{\frac{1}{2}}$-adjoint operators. Furthermore, several properties of $A$-normal operators have been established.
	\end{abstract}
	
      \begin{amscode} 47A10; 47A30; 47A80; 46C05.\end{amscode}
    
    \section{Introduction}
  Let $H$ be a Hilbert space over the complex field. The algebra of all linear bounded operators on $H$ is denoted by $B(H)$. We consider $A$ as a non-zero positive semidefinite operator. The sesquilinear form ${\langle x , y \rangle}_{A}$  is defined from $H \times H$ to $\mathbb C$ by ${\langle x , y \rangle}_{A} = \langle Ax , y \rangle$. The associated seminorm ${\| .\|}_{A}$ is derived from  ${\langle . ,. \rangle}_{A}$, expressed as  ${\| x \|}_{A} = {\langle x , x \rangle}_{A}^{\frac{1}{2}}$. Notably, $\| x\|_{A} = 0$ if and only if $x\in N(A)$, where $N(A)$ denotes the null space of $A$.  The subspace ${\mathcal{M}}^{\perp_{A}} = \{x : \langle Ax, y \rangle = 0 \text{ for all } y \in \mathcal{M}\}$ is termed the $A$-orthogonal companion of $\mathcal{M}$. $R(T)$ denotes the range of an operator $T$, and its closure is denoted by $\overline{R(T)}$.

This paper delves into the exploration of $A$-approximate point spectrum for $A$-bounded operators acting on a complex semi-Hilbertian space $H$, elucidating several characteristics of $A$-normal operators. It is very much obvious that  we can express $\|Sx\| = \|x\|_{|S|^{2}}$ for any bounded operator $S$ in $H$. This concept allows the utilization of properties of semi-Hilbertian spaces in the analysis of bounded operators in $H$. The organizational structure of the paper unfolds in three main sections. Section 2 delves into various characterizations of the $A$-approximate point spectrum, while Section 3 investigates properties of the $A$-approximate point spectrum for the tensor product of two $A^{\frac{1}{2}}$-adjoint operators. The final section establishes characterizations of $A$-normal operators, providing a comprehensive framework for understanding the spectral properties of $A$-bounded operators in complex semi-Hilbertian spaces.
 \begin{definition}\cite{MR2388631} In the context of Hilbert spaces, an operator  $T$ belonging to $B(H)$ is termed an $A$-bounded operator if there exists a constant $d > 0$ such that for all $\xi$ in the closure of the range of $A$, the inequality $\|T \xi\|_{A} \leq d  \|\xi\|_{A}$ holds. The $A$-norm of $T$ is defined as
 \begin{equation*}
     {\|T\|}_{A} = \sup_{\xi \in \overline {R(A)}\setminus \{0\}}\frac{{\|T \xi\|}_{A}}{{\| \xi\|}_{A}} < \infty.
 \end{equation*}
 Alternatively, this norm can be expressed as
 \begin{align*}
 {\|T\|}_{A} = \sup \{ | {\langle T \xi, \eta \rangle}_{A} | : \xi, \eta \in H , {\| \xi\|}_{A} \leq 1, {\| \eta\|}_{A} \leq 1\}.
 \end{align*}
 The set of all $A$-bounded operators is denoted by $B^{A}(H)= \{T \in B(H) : {\|T\|}_{A} < \infty\}$.

 \end{definition}
 \begin{theorem} \cite{MR0203464}\label{thm 1}
 Let $E_{1}, E_{2} \in B(H)$. The following conditions are equivalent:
 \begin{enumerate}
     \item $R(E_{2}) \subset R(E_{1})$.
     \item There exists a positive number $\mu $ such that $E_{2}E_{2}^{*} \leq \mu E_{1}E_{1}^{*}$.
     \item There exists $C\in B(H)$ such that $E_{1}C = E_{2}$.
 \end{enumerate}
 If at least one of these conditions holds, a unique operator $D \in B(H)$ exists such that $E_{1}D = E_{2}$ with $R(D) \subset \overline {R(E_{1}^{*})}$. Furthermore, $N(D) = N(E_{2})$, and $D$ is referred to as the reduced solution of the equation $E_{1}X = E_{2}$.
 \end{theorem}
 \begin{definition}\cite{MR2442900}
Consider a bounded linear operator $V$ as an $A$-adjoint of $T\in B(H)$ if, $\text{ for all }   \xi, \eta \in H$, it satisfies the following condition
 \begin{align*}
 {\langle T\xi , \eta \rangle}_{A} = {\langle \xi, V\eta \rangle}_{A}.
 \end{align*}
 \end{definition}
 Building on Theorem (\ref{thm 1}), the existence of an $A$-adjoint operator for $T\in B(H)$ is established if and only if $R(T^{*}A) \subset R(A)$. We define $B_{A}(H)$ as the subalgebra of $B(H)$ comprising operators with $A$-adjoint operators, that is,
 \begin{align*}
 B_{A}(H) = \{T \in B(H) : R(T^{*}A) \subset R(A)\}.
 \end{align*}
 Similarly, we introduce ${B_{A^{\frac{1}{2}}}}(H) = \{ T \in B(H) : R(T^{*}A^{\frac{1}{2}}) \subset R(A^{\frac{1}{2}})\}$. As per Theorem (\ref{thm 1}), it is evident that 
 \begin{align*}
 B_{A^{\frac{1}{2}}}(H) = \{ T \in B(H) : \exists ~d>0, ~{\|T\xi\|}_{A} \leq d {\|\xi\|}_{A},  \text { for all } \xi \in H\}.
 \end{align*}
In \cite{MR2590353}, the inclusion relations  are established as $B_{A}(H) \subset B_{A^{\frac{1}{2}}}(H) \subset B^{A}(H) \subset B(H)$.
 
For $T\in B_{A}(H)$, there exists a unique A-adjoint denoted as $T^{\sharp}$, satisfying $T^{*}A = A T^{\sharp}$ with $R(T^{\sharp}) \subset \overline {R(A)}$. The reduced solution $T^{\sharp}$ of $T^{*}A= AX$ implies $N(T^{\sharp}) = N(T^{*}A)$, and $T^{\sharp} = A^{\dagger}T^{*}A$, where $A^{\dagger}$ represents the Moore-Penrose of $A$ in the domain $ D(A^{\dagger}) =R(A) \oplus R(A)^{\perp}$.

For $T \in B_{A}(H)$ termed $A$-selfadjoint when $T^{*}A = AT$. In a similar way, if $T \in B_{A^{\frac{1}{2}}}(H)$, there exists a unique reduced solution $T^{\diamond}$ such that $T^{*}A^{\frac{1}{2}} = A^{\frac{1}{2}} T^{\diamond}$ with $R(T^{\diamond}) \subset \overline{R(A^{\frac{1}{2}})}$.
\begin{proposition}\cite{MR2442900}\label{ pro 2.5}
Let $T\in B(H)$. Then the following statements are equivalent:
\begin{enumerate}
    \item $T \in B^{A}(H)$.
    \item $A^{\frac{1}{2}}T(A^{\frac{1}{2}})^{\dagger}$ is bounded.
    \item $R(A^{\frac{1}{2}}T^{*}A^{\frac{1}{2}}) \subset R(A)$.
\end{enumerate}
Furthermore, one of these conditions ensures that ${\|T\|}_{A} = \|A^{\frac{1}{2}}T(A^{\frac{1}{2}})^{\dagger}\| = \|(A^{\frac{1}{2}})^{\dagger} T^{*}A^{\frac{1}{2}}\|$.
\end{proposition}
Here, we delve into several noteworthy properties associated with the operator $T^{\sharp}$ extensively explored in the literature \cite{MR2442900, MR2388631}. For the sake of completeness, we present and elaborate on these properties. The notation $P$ replacing $P_{\overline{R(A)}}$ throughout,  where $P_{\overline{R(A)}}$ signifies the orthogonal projection onto the closure of the range of $A$.
\begin{proposition}\cite{MR2442900, MR2388631}
    Let $T \in B_{A}(H)$. Then the following statements are true:
    \begin{enumerate}
        \item  For any $c > 0$, it is established that $(A^{c})^{\sharp} = A^{c}$.
        \item  In a situation where $AT = TA$, the relationship  $T^{\sharp} = PT^{*}$ holds.
        \item If $AT = T^{*}A$, then $(A^{\frac{1}{2}})^{\dagger}T^{*}A^{\frac{1}{2}}$ manifests as a  selfadjoint operator.
        \item If $W \in B_{A}(H)$, then $TW$ is an element of $B_{A}(H)$ and the relation $(TW)^{\sharp}= W^{\sharp}T^{\sharp}$ holds good.
        \item The operator $T^{\sharp}$ is in $B_{A}(H)$ and $(T^{\sharp})^{\sharp} = PTP$. Notably, $((T^{\sharp})^{\sharp})^{\sharp} = T^{\sharp}$.
        \item ${\|T\|}_{A} = {\|T^{\sharp}\|}_{A} = {{\|T^{\sharp}T\|}_{A}}^{\frac{1}{2}}$.
\end{enumerate}
\end{proposition}
\begin{definition}\cite{MR4133606}
Let $T \in B(H)$. The $A$-numerical range, denoted by $W_{A}(T)$, is defined as follows:
\begin{align*}
W_{A}(T) = \{\langle Tx, x\rangle_{A} : x \in H, \|x\|_{A} = 1\}.
\end{align*}
Furthermore, the $A$-numerical radius, denoted by $w_{A}(T)$, is defined as the supremum of absolute values: $w_{A}(T) = \sup \{|\langle Tx, x\rangle|_{A} : x\in H, \|x\|_{A}= 1\}$.
\end{definition}
\begin{theorem}\cite{MR3834203}\label{thm 2.14}
For any $T \in B(H)$, the $A$-numerical range $W_{A}(T)$ forms a convex subset of the complex plane $\mathbb C$.
\end{theorem}
 
\begin{definition}\label{def 2.11}
 Let $T \in B(H)$. A scalar $\lambda \in \mathbb {C}$ is said to be  an $A$-point spectral value of $T$ if there exists $x\in H$ with $\|x\|_{A} \neq 0$ such that $x \in N(A^{\frac{1}{2}}(\lambda -T))$.
 The set of all $A$-point spectral values of $T$ is termed the $A$-point spectrum of $T$, denoted by $\sigma_{A_{p}}(T)$.
 \end{definition}
 \begin{definition}\label{def 2.12}
 Let $T \in B(H)$. A scalar $\lambda \in \mathbb {C}$ is said to be an $A$-approximate point spectral value of $T$ if there exists a sequence $\{x_{n}\}$ in $H$, where $\|x_{n}\|_{A}= 1$, such that $\|(T- \lambda) x_{n}\|_{A} \to 0 \ \  \text{as} \  \ n \to \infty.$ 
\end{definition} 
The set of all $A$-approximate point spectral values of $T$ is called the $A$-approximate point spectrum of $T$, denoted by $\sigma_{A_{app}}(T)$. 

Notably, if $T \in B_{A^{\frac{1}{2}}}(H)$,  $\lambda \in \sigma_{A_{app}}(T)$ can equivalently be characterized by the existence of a sequence $\{z_{n}\} \text{ in } \overline{R(A)}$ with $\|z_{n}\|_{A} = 1$ such that $\|(T -\lambda)z_{n}\|_{A} \to 0$ as $n \to \infty$. Additionally, it is evident that  $\sigma_{A_{p}}(T) \subset \sigma_{A_{app}}(T)$, for all $T\in B(H)$.

The notions of $A$-point and $A$-approximate point spectral values are studied in \cite{arxiv 007} for $A^{\frac{1}{2}}$-adjoint operators. Our mentioned definitions are generalized versions of any bounded operator.

 \begin{definition}\cite{MR4362420}
  A non-zero operator $T\in B_{A^{\frac{1}{2}}}(H)$ is termed as $A$-invertible in $ B_{A^{\frac{1}{2}}}(H)$ if there exists a non-zero operator $S\in B_{A^{\frac{1}{2}}}(H)$ such that $ATS = AST =A$. The operator $S$ is referred to as an $A$-inverse in $ B_{A^{\frac{1}{2}}}(H) $.  
  
  Similarly, a non-zero operator $T \in  B_{A}(H)$ is said to be $A$-invertible in $ B_{A}(H)$ if there exists a non-zero operator $S \in  B_{A}(H)$ such that $ATS = AST =A$. Here, $S$ is called an $A$-inverse in $B_{A}(H)$.
 \end{definition}
 \begin{definition}\cite{MR4362420}
 Consider $T \in B_{A^{\frac{1}{2}}}(H)$.
 \begin{enumerate}
     \item The $A$-resolvent set of $T$, denoted by $\rho_{A}(T)$, is defined as $\rho_{A}(T) = \{ \lambda \in \mathbb {C}: (\lambda - T)$ is $A$-invertible in $B_{A^{\frac{1}{2}}}(H)\}$.
     \item The $A$-spectrum of $T$ is denoted by  $\sigma_{A}(T) = \mathbb{C} \setminus \rho_{A}(T)$.
 \end{enumerate}
 \end{definition}
  \begin{remark}\cite{MR4362420}
 For $T \in B_{A^{\frac{1}{2}}}(H)$ being an $A$-invertible operator in $B_{A^{\frac{1}{2}}}(H)$ with an $A$-inverse $S \in B_{A^{\frac{1}{2}}}(H)$, the following statements are equivalent:
 \begin{enumerate}
   \item $ATS = AST =A$.
   \item $PTS = PST =P$.
   \item $A^{\frac{1}{2}}TS = A^{\frac{1}{2}}ST = A^{\frac{1}{2}}$.
 \end{enumerate}
 \end{remark}
 \begin{remark}\label{remark 1}
 Suppose $T, S$ are operators in $B_{A}(H)$. The condition for $T$ to be an $A$-invertible operator in $B_{A}(H)$ with an $A$-inverse $S \text { in } B_{A}(H)$ is equivalent to $T^{\sharp}$ possessing same property with an $A$-inverse $S^{\sharp} \text { in } B_{A}(H)$.
 \end{remark}
\begin{definition}\cite{MR4133606}
For $T \in B_{A^{\frac{1}{2}}}(H)$, the $A$-spectral radius of $T$ is defined as
\begin{align*}
r_{A}(T) = \lim_{n \to \infty} (\|T^{n}\|_{A})^{\frac{1}{n}}.
\end{align*}
\end{definition}
\begin{theorem}\cite{MR4133606}\label{thm 2.15}
If $T \in B_{A^{\frac{1}{2}}}(H)$, then $r_{A}(T)$ satisfies $r_{A}(T) \leq w_{A}(T) \leq \|T\|_{A}$.
\end{theorem}
 \begin{theorem}\cite{hb}\label{thm 2}
 If $T \in B_{A}(H)$  is $A$-invertible in $B_{A^{\frac{1}{2}}}(H)$, any $A$-inverse of $T$ in $B_{A^{\frac{1}{2}}}(H)$ also belongs to $B_{A}(H)$.
 \end{theorem}
\begin{definition}\cite{MR2590353}
Let $T\in B(H)$. The $A$-reduced minimum modulus of $T$ is given by 
    \begin{align}
    {\gamma_{A}}(T) = \inf \{\|T\xi\|_{A}: \xi \in {N(A^{\frac{1}{2}}T)^{\perp_{A}} , \|\xi\|_{A}=1}\}.
    \end{align}
For $T \in B_{A}(H)$, ${\gamma_{A}}(T) =  \inf \{\|T\xi\|_{A}: \xi \in \overline{R(T^{\sharp}T)} , \|\xi\|_{A}=1\}$.
\end{definition}
\begin{proposition}\cite{MR2590353}
Consider $T \in B_{A}(H)$ and a solution $E$ of the equation $A^{\frac{1}{2}}X = T^{*}A^{\frac{1}{2}}$. If $A^{\frac{1}{2}} \overline{R(T^{\sharp}T)} \subset \overline {R(E)}$, then $\gamma_{A}(T) = \gamma(E)$.
\end{proposition}
\begin{proposition}\cite{MR2590353}\label{pro 0.1}
Let $T \in B_{A}(H)$. Then
\begin{enumerate}
 \item $\gamma_{A}(T) = \gamma(T^{\diamond})$.
 \item $\gamma_{A}(T) = \gamma_{A}(T^{\sharp})$.
\end{enumerate}
\end{proposition}
 Now let us consider the Hilbert space $\mathbf{R}(A^{\frac{1}{2}}) = (R(A^{\frac{1}{2}}), (.,.))$  with the inner product $(A^{\frac{1}{2}} x, A^{\frac{1}{2}} y) = \langle Px, Py\rangle,  \text {for all } x, y \in H$. Then $\|A^{\frac{1}{2}} x\|_{\mathbf{R}(A^{\frac{1}{2}})} = \|Px\|, \text { for all } x \in H$. We define two operators: $W_{A} : H \mapsto \mathbf{R}(A^{\frac{1}{2}})$  by $W_{A}(x) = Ax, \text { for all } x\in H$ and $Z_{A}: H \mapsto \mathbf{R}(A^{\frac{1}{2}})$ by $Z_{A}(x) = A^{\frac{1}{2}}(x), \text{ for all } x\in H$. Several properties of the operators $W_{A} \text{ and } Z_{A}$ have been studied extensively in the paper \cite{MR2590353}.
 
 The subsequent outcome delineates the connection between $A$-bounded operators in a semi-Hilbertian space $H$ and operators in $B(\mathbf{R}(A^{\frac{1}{2}}))$.
\begin{proposition}\cite{MR2590353}
Let $T \in B(H)$. The operator $T\in B_{A^{\frac{1}{2}}}(H)$ if and if there exists an unique $\tilde{T} \in B(\mathbf{R}(A^{\frac{1}{2}}) )$ such that $\tilde{T}W_{A} = W_{A}T$.
\end{proposition}

\section{ Some characterizations of the $A$-approximate point spectrum}
The study conducted by Hamadi Baklouti and Sirine Namouri in 2021 \cite{MR4362420} delves into the spectral analysis of bounded operators on semi-Hilbertian spaces. In this context, we present various characterizations of the $A$-approximate point spectrum of $T \in B_{A^{\frac{1}{2}}}(H)$. The following results contribute to a comprehensive understanding of the properties and relationships within this mathematical framework.
 \begin{theorem}\label{thm 0.1}
     Let $T\in B_{A}(H)$ be $A$-invertible in $B_{A}(H)$. Then ${\gamma_{A}}(T) = {\gamma_{A}}(T^{\sharp}) > 0$.
 \end{theorem}
 \begin{proof}  $A$-invertibility of $T\in B_{A}(H)$ in $B_{A}(H)$ confirms that there exists $S\in B_{A}(H)$ such that  $AST = ATS =A$. Moreover, $\|S\|_{A} \neq 0$ because $\|S\|_{A} = 0 \text{ if and only if } AS= 0$ which implies $A= AST = 0$. 
     \begin{align*}
         \|x\|_{A}= \|STx\|_{A} \leq \|S\|_{A}\|Tx\|_{A} , \  \text{ for all }\ x\in H.
     \end{align*}
     Therefore, $0 < \frac{1}{\|S\|_{A}} \leq \gamma_{A}(T) = \gamma_{A}(T^{\sharp})$.
 \end{proof}
 \begin{remark}
 Let $T \in B_{A^{\frac{1}{2}}}(H)$ be $A$-invertible in $B_{A^{\frac{1}{2}}}(H)$. Then $0 < \frac{1}{\|S\|_{A}} \leq \gamma_{A}(T)$. Moreover, $\frac{1}{\|S\|_{A}} = \gamma_{A}(T)$ when $\|S\|_{A}\|T\|_{A} = 1$.
 \end{remark}
\begin{corollary}
 Let $T \in B_{A}(H)$ be $A$-invertible in $B_{A}(H)$. Then $R(T^{\diamond})$ is closed.
 \end{corollary}
 \begin{proof}
Proposition \ref{pro 0.1} and Theorem \ref{thm 0.1} guarantee that $\gamma(T^{\diamond}) > 0$. Therefore $R(T^{\diamond})$ is closed.
 \end{proof}
 
    The following Theorem is similarly proved as Proposition 2.7\cite{arxiv 007}, but we consider $T \in B(H)$.

 \begin{theorem}\label{thm 3.32}
Let $T \in B(H)$. Then $\sigma_{A_{app}}(T) \subset \overline{W_{A}(T)}$.
\end{theorem}
\begin{proof}
Let us consider $\lambda \in \sigma_{A_ {app}}(T)$. Then there exists a sequence $\{x_{n}\}$ in the Hilbert space $H$, where $\|x_{n}\|_{A}= 1$, such that $\|(\lambda -T)x_{n}\|_{A} \to 0$ as $n \to \infty$. Then
\begin{align}\label{ineq 3.32}
|\langle (T- \lambda)x_{n}, x_{n}\rangle_{A}|  \leq \|(T- \lambda)x_{n}\|_{A} \|x_{n}\|_{A}.
\end{align}
The right-hand side of the above inequality (\ref{ineq 3.32}) goes to $0$ as $n \to \infty$. Thus $\lambda \in \overline{W_{A}(T)}$. Therefore,  $\sigma_{A_{app}}(T) \subset \overline{W_{A}(T)}$.
\end{proof}

\begin{theorem}\label{thm 4}
Let $T \in B(H)$. Then $\sigma_{A_{app}}(T)$ is closed.
\end{theorem}
\begin{proof}
Let us consider an arbitrary element $\lambda \in \overline{\sigma_{A_{app}}(T)}$. Then there exists a sequence $\{\lambda_{n}\}$ in $\sigma_{A_{app}}(T)$ such that $\lambda_{n} \to \lambda$ as $n \to \infty$. If $\lambda \notin \sigma_{A_{app}}(T)$, then there exists a constant $c>0$ such that $\|(T-\lambda)x\|_{A} \geq c \|x\|_{A}, \text { for all } x\in H$. Moreover, we can get a natural number $m\in \mathbb N$ such that $|\lambda_{n}- \lambda| < \frac{c}{2}, \text{ for all } n\geq m$. Thus,
\begin{equation}
\|(T-\lambda_{m})x\|_{A} = \|(T- \lambda)x + (\lambda - \lambda_{m})x\|_{A} \geq \frac{c}{2}\|x\|_{A}, \text { for all } x\in H.
\end{equation}
This is a contradiction because $\lambda_{m} \in \sigma_{A_{app}}(T)$.  Therefore, $\sigma_{A_{app}}(T)$ is closed.
\end{proof}
\begin{theorem}\label{thm 5}
Let $T \in  B_{A^{\frac{1}{2}}}(H)$. Then ${\sigma_{A_{app}}(T)} \subset \sigma_{A}(T)$.
\end{theorem}
\begin{proof}
Let us consider an arbitrary element $\lambda \in \sigma_{A_{app}}(T)$ but $\lambda \notin \sigma_{A}(T)$. Then there exists an $A$-inverse $S_{\lambda} \in B_{A^{\frac{1}{2}}}(H)$ of $(\lambda -T) \in B_{A^{\frac{1}{2}}}(H)$ such that $A^{\frac{1}{2}}(\lambda - T)S_{\lambda} = A^{\frac{1}{2}}S_{\lambda}(\lambda - T) = A^{\frac{1}{2}}$. Additionally, we obtain a sequence $\{x_{n}\}$ in $H$, satisfying $\|x_{n}\|_{A} =1$, such that the expression $\|(\lambda -T)x_{n}\|_{A} \to 0$ as $n \to \infty$. Also, we get
\begin{equation}\label{ineq 2}
\|x_{n}\|_{A} = 1 = \|A^{\frac{1}{2}} S_{\lambda} (\lambda -T)x_{n}\| = \|({S_{\lambda}}^{\diamond})^{*} A^{\frac{1}{2}}(\lambda -T)x_{n}\| \leq \|({S_{\lambda}}^{\diamond})^{*}\| \| (\lambda -T)x_{n}\|_{A}.
\end{equation}
The left side of the above inequality (\ref{ineq 2}) is $1$, but the right-hand side goes to $0$ as $n \to \infty$, which is a contradiction. Therefore, $\sigma_{A_{app}}(T) \subset \sigma_{A}(T)$.
\end{proof}
\begin{corollary}
Let $T \in B_{A^{\frac{1}{2}}}(H)$. If the conditions $w_{A}(T) = \|T\|_{A}$ and $T A^{\frac{1}{2}} = A^{\frac{1}{2}}T$ hold, then $r_{A}(T) = \|T\|_{A}$.
\end{corollary}
\begin{proof}
Without loss of generality, we can assume that  $w_{A}(T) = \|T\|_{A} = 1$. Then there exists $\lambda \in \mathbb C$ with $|\lambda| =1$ such that $\langle Tx_{n}, x_{n}\rangle_{A} \to \lambda$ as $n \to \infty$, for some sequence $\{x_{n}\}$ in $H$ with $\|x_{n}\|_{A}= 1$. By Cauchy-Schwarz inequality, $|\langle Tx_{n}, x_{n}\rangle_{A}| \leq \|Tx_{n}\|_{A} \leq 1$. So, $\|Tx_{n}\|_{A} \to 1$ as $n \to \infty$.
\begin{align}\label{eq 20}
\|(T - \lambda)x_{n}\|_{A}^{2} = \|Tx_{n}\|_{A}^{2} - \lambda \langle Ax_{n}, Tx_{n}\rangle - \overline{\lambda} \langle ATx_{n}, x_{n}\rangle + |\lambda|^{2} \langle Ax_{n}, x_{n}\rangle.
\end{align}
The right-hand side of the equality (\ref{eq 20}) goes to $0$ as $n \to \infty$. Thus, $\lambda \in \sigma_{A_{app}}(T) \subset \sigma_{A}(T)$. By Theorem (3.11) \cite{arxiv 001} and Theorem (\ref{thm 2.15}), we get $\|T\|_{A} = 1 \leq r_{A}(T) \leq w_{A}(T) = \|T\|_{A}$. Therefore, $r_{A}(T) = \|T\|_{A}$.
\end{proof}
\begin{theorem}
Let $T\in B_{A^{\frac{1}{2}}}(H)$. Then $\sigma_{A_{app}}(T) = \sigma_{A_{app}}(PT) = \sigma_{A_{app}}(TP)$.
\end{theorem}
\begin{proof}
The scalar $\lambda$ belongs to $\sigma_{A_{app}}(T)$ if and only if we obtain a sequence $\{x_{n}\}$ in $\overline{R(A)}$, where $\|x_{n}\|_{A} = 1$, such that $\|(\lambda -T) x_{n}\|_{A} \to 0$ as $n \to \infty$. Equivalently this conditions holds if there exists a sequence $\{x_{n}\}$ in $\overline{R(A)}$ with $\|x_{n}\|_{A} = 1$ such that $\| A^{\frac{1}{2}} (\lambda - PT) x_{n}\| = \|(\lambda - PT) x_{n}\|_{A}\to 0$ as $n \to \infty$. Consequently it follows that $\sigma_{A_{app}}(T) = \sigma_{A_{app}}(PT)$.
Similarly, the same reasoning can be applied to demonstrate the second equality.
\end{proof}
\begin{theorem}
Let $T \in B_{A^{\frac{1}{2}}}(H)$ be  $A$-invertible in $B_{A^{\frac{1}{2}}}(H)$ with $A$-inverse $S \in B_{A^{\frac{1}{2}}}(H)$. Then $\sigma_{A_{app}}(S) = \{\lambda \in \mathbb C : \frac{1}{\lambda} \in \sigma_{A_{app}}(T)\}$.
\end{theorem}
\begin{proof}
It is obvious that $0 \in \rho_{A}(T)$ and $0 \in \rho_{A}(S)$. So, $0 \notin \sigma_{A_{app}}(T)$ and $0 \notin \sigma_{A_{app}}(S)$. Let $\frac{1}{\lambda} \in \sigma_{A_{app}}(T)$. Then there exists a sequence $\{x_{n}\}$ in $H$, where $\|x_{n}\|_{A} =1$, such that 
\begin{align}\label{eq 1}
\Big\|\Big(\frac{1}{\lambda} - T\Big)x_{n}\Big\|_{A} = \frac{\| (S - \lambda)Tx_{n}\|_{A}}{|\lambda|}.
\end{align}
The equality (\ref{eq 1}) goes to $0$ as $n \to \infty$. Since $0 \notin \sigma_{A_{app}}(T)$, we get a constant $a > 0$ such that $\|Tx\|_{A} \geq a\|x\|_{A}$. So, $\frac{1}{\|Tx_{n}\|_{A}} \leq \frac{1}{a},  \text { for all } n\in \mathbb N$. Now consider, $z_{n} = \frac{Tx_{n}}{\|Tx_{n}\|_{A}}, \text { for all } n\in \mathbb N$.
\begin{align}\label{ineq 6}
\|(\lambda - S)z_{n}\|_{A} = \frac{\|(\lambda - S)Tx_{n}\|_{A}}{\|Tx_{n}\|_{A}} \leq \frac{1}{a} \|(\lambda - S)Tx_{n}\|_{A}.
\end{align}
The right-hand side of the inequality (\ref{ineq 6}) goes to $0$ as $n \to \infty$.
Thus, $\lambda \in \sigma_{A_{app}}(S)$.
  
  Now, we will prove the converse part. Let us consider $\mu \in\sigma_{A_{app}}(S)$. Then we obtain a sequence $\{z_{n}\}$ in $H$, where $\|z_{n}\|_{A} =1$, such that $\|(\mu - S)z_{n}\|_{A} = |\mu|{\|(\frac{1}{\mu} -T)Sz_{n}\|}_{A} \to 0$ as $n \to \infty$. Since $0 \notin \sigma_{A_{app}}(S)$, we get a positive real number $d >0$ such that $\|Sz\|_{A} \geq d\|z\|_{A}$. So, $\frac{1}{d} \geq \frac{1}{\|Sz_{n}\|_{A}}, \text { for all } n\in \mathbb N$. Consider, $w_{n} = \frac{Sz_{n}}{\|Sz_{n}\|_{A}} , \text { for all } n\in \mathbb N$.
  \begin{align}\label{ineq 7}
  \Big\|\Big(\frac{1}{\mu} -T\Big)w_{n}\Big\|_{A} = \frac{\|(\frac{1}{\mu} - T)Sz_{n}\|_{A}}{\|Sz_{n}\|_{A}} \leq {\frac{1}{d}} \Big\|\Big(\frac{1}{\mu} -T\Big)Sz_{n}\Big\|_{A}.
  \end{align}
  The right-hand side of the inequality (\ref{ineq 7}) goes to $0$ as $n \to \infty$. Hence, $\frac{1}{\mu} \in \sigma_{A_{app}}(T)$. Therefore, $\sigma_{A_{app}}(S) = \{\lambda \in \mathbb C : \frac{1}{\lambda} \in \sigma_{A_{app}}(T)\}$.
\end{proof}
\begin{remark}
Let $W \in B_{A^{\frac{1}{2}}}(H)$ be $A$-invertible  with an $A$-inverse $V \text { in } B_{A^{\frac{1}{2}}}(H)$. Let $T$ be in $B_{A^{\frac{1}{2}}}(H)$. Then, $\sigma_{A_{app}}(WVT) = \sigma_{A_{app}}(VWT) = \sigma_{A_{app}}(T)$.
\end{remark}
\begin{theorem}
Let $T\in B_{A^{\frac{1}{2}}}(H)$. If $\lambda \in \sigma_{A_{app}}(T)$, then ${\lambda}^{n} \in \sigma_{A_{app}}(T^{n}), \text{ for all } n\in \mathbb N$.
\end{theorem}
\begin{proof}
For $\lambda \in \sigma_{A_{app}}(T)$, we obtain a sequence $\{x_{n}\}$ in $H$, satisfying $\|x_{n}\|_{A} =1$, such that $\|(\lambda -T)x_{n}\|_{A} \to 0$ as $n \to \infty$. Then,
\begin{align}\label{eq 0.2}
\|({\lambda}^{n} - T^{n})x_{n}\|_{A} \leq \| {\lambda}^{n-1} + \dots + {\lambda}T^{n-2} + T^{n-1}\|_{A}  \|(\lambda -T)x_{n}\|_{A}.
\end{align}
The right-hand of the inequality (\ref{eq 0.2}) goes to $0$ as $n \to \infty$. Therefore, ${\lambda}^{n} \in \sigma_{A_{app}}(T^{n}),  \text { for all } n\in \mathbb N$.
\end{proof}

\begin{theorem}
Let $T \in B_{A^{\frac{1}{2}}}(H)$. If $\lambda \in W_{A}(T)$ with $|\lambda| = \|T\|_{A}$, then $\lambda \in \sigma_{A_{p}}(T)$.
\end{theorem}
\begin{proof}
For $\lambda \in W_{A}(T)$, there is an element  $x \in H$, where $\|x\|_{A} = 1$, such that $\lambda = \langle Tx, x\rangle_{A}$. So, $\|T\|_{A} = |\lambda| = |\langle A^{\frac{1}{2}}Tx, A^{\frac{1}{2}}x \rangle| \leq \|Tx\|_{A} \|x\|_{A} \leq \|T\|_{A}$. Thus, $|\langle A^{\frac{1}{2}}Tx, A^{\frac{1}{2}}x \rangle| = \|A^{\frac{1}{2}}Tx\| \|A^{\frac{1}{2}}x\|$. There exists $\mu \in \mathbb C$ with $A^{\frac{1}{2}}Tx = \mu A^{\frac{1}{2}}x$ which implies $ATx = \mu Ax$. Now, $\lambda = \langle \mu Ax, x\rangle = \mu$ implies $\|(T- \lambda)x\|_{A} = 0$. Therefore, $\lambda \in \sigma_{A_{p}}(T)$. 
\end{proof}
\begin{theorem}\cite{hb}\label{thm 3.37}
Under the assumption that $A$ has a closed range and $T \in  B_{A^{\frac{1}{2}}}(H)$, the $A$-invertibility of $T$ in $B_{A^{\frac{1}{2}}}(H)$ is established if and only if $\tilde T$ is invertible in $B(\mathbf{R}(A^{\frac{1}{2}}))$. Specifically, this equivalence extends to the equality as $\sigma_{A}(T) = \sigma(\tilde{T}), \text { for all } T\in B_{A^{\frac{1}{2}}}(H)$. Here $\sigma(\tilde{T})$ denotes the spectrum of $\tilde{T}$. 
\end{theorem}
\begin{corollary}
In a finite-dimensional Hilbert space $H$, $\sigma_{A}(T)= \sigma_{A_{app}}(T) = \sigma_{A_{p}}(T)$ for all $T \in B_{A^{\frac{1}{2}}}(H)$.
\end{corollary}
\begin{proof}
This result follows from Theorem 2.13 \cite{arxiv 007} and Proposition 2.5 \cite{arxiv 007}.
\end{proof}
The closed range of the operator $A$ is used in the proof of Theorem 5.1\cite{MR4362420}, though it is not mentioned in the statement of Theorem 5.1\cite{MR4362420}. Moreover, the correct version of Theorem 5.1\cite{MR4362420} has been established in Theorem 3.1\cite{arxiv 001}. In literature, Theorem 5.1\cite{MR4362420} has been used to prove several other results including Theorem 5.3\cite{MR4362420} and Theorem 2.11 \cite{arxiv 007}. The modified version of Theorem 5.3\cite{MR4362420} is proved in \cite{arxiv 001} and it is stated in the following Lemma.
\begin{lemma}\label{lemma 1} \cite{arxiv 001}
Let $T\in B_{A^{\frac{1}{2}}}(H)$ be invertible in $B_{A^{\frac{1}{2}}}(H)$, possessing an $A$-inverse denoted by $S$. Let $T^{\prime}$ be another operator in $ B_{A^{\frac{1}{2}}}(H)$ such that $\|T^{\prime}S\|_{A} < 1 \text{ and } \|(T^{\prime}S)^{\diamond}\|_{A} < 1 $. Then $T + T^{\prime}$ is $A$-invertible in $B_{A^{\frac{1}{2}}}(H)$.
\end{lemma}
\begin{theorem}\label{cor 1}
If $T$ belongs to the space $B_{A^{\frac{1}{2}}}(H)$, then both the $A$-spectrum $\sigma_{A}(T)$ and $A$-approximate point spectrum $\sigma_{A_{app}}(T)$ are compact sets within the complex plane $\mathbb C$. 
\end{theorem}
\begin{proof}
First, we will show that $\rho_{A}(T)$ is open. Let us consider an arbitrary element $\lambda_{0} \in \rho_{A}(T)$. Then there exists an $A$-inverse of $(\lambda_{0} -T)$ in $B_{A^{\frac{1}{2}}}(H)$, say $S_{\lambda_{0}}$, such that $P(\lambda_{0} -T) S_{\lambda_{0}} = PS_{\lambda_{0}}(\lambda_{0} - T) = P$. It is obvious from Lemma \ref{lemma 1} that $\lambda -T$ is invertible on $B_{A^{\frac{1}{2}}}(H) \text { for all } \lambda \in \{\mu : |\mu - \lambda_{0}| < \frac{1}{max\{\|S_{\lambda_{0}}\|_{A}, \|S^{\diamond}_{\lambda_{0}}\|_{A}\}}\}$. Thus, $\rho_{A}(T)$ is open.
    Now Corollary 3.15 \cite{arxiv 001} confirms that $\sigma_{A}(T)$ is compact. Furthermore, Theorems \ref{thm 4} and  \ref{thm 5} conclude that $\sigma_{A_{app}}(T)$ is compact.
\end{proof}
\begin{theorem}\label{thm 6}
Let $T\in B_{A^{\frac{1}{2}}}(H)$. Then the conjugate set $(\sigma_{A}(T^{\diamond}))^{*}$ is equal to the $A$-spectrum set $\sigma_{A}(T)$, where $(\sigma_{A}(T^{\diamond}))^{*}$ denotes the conjugate set of $\sigma_{A}(T^{\diamond})$.
\end{theorem}
\begin{proof}
 $\lambda \in \rho_{A}(T)$ if and only if there exists $S_{\lambda} \in B_{A^{\frac{1}{2}}}(H)$ such that $P(\lambda -T) S_{\lambda} = PS_{\lambda}(\lambda -T) = P$ if and only if $P(S_{\lambda})^{\diamond}(\overline{\lambda}P- T^{\diamond}) = P (\overline{\lambda}P -T^{\diamond})(S_{\lambda})^{\diamond} = P$ if and only if $P(S_{\lambda})^{\diamond}(\overline{\lambda}- T^{\diamond}) = P (\overline{\lambda}- T^{\diamond})(S_{\lambda})^{\diamond} =P$ if and only if $\overline{\lambda} \in \rho_{A}(T^{\diamond})$. Therefore, $(\sigma_{A}(T^{\diamond}))^{*} = \sigma_{A}(T)$.
\end{proof}
\begin{corollary} Let $T \in B_{A}(H)$.  Then  $(\sigma_{A}(T^{\sharp}))^{*} = \sigma_{A}(T)$.
\end{corollary}
\begin{theorem}\label{thm 7}
Let $T\in B_{A^{\frac{1}{2}}}(H)$. Then $\sigma_{A_{app}}(T)= \sigma_{A_{app}}((T^{\diamond})^{\diamond})$.
\end{theorem}
\begin{proof}
For any $\lambda \in \sigma_{A_{app}}(T)$, there exists $\{x_{n}\} \in \overline{R(A)}$ with $\|x_{n}\|_{A} = 1$ such that $\|(\lambda - T)x_{n}\|_{A} \to 0$ as $n \to \infty$. The condition is equivalent to the existence of a sequence $\{x_{n}\} \in \overline{R(A)}$ with $\|x_{n}\|_{A} = 1$ such that $\|(\lambda P - PTP)x_{n}\|_{A} \to 0$ as $n \to \infty$. Moreover, it is also equivalent to the existence of a sequence $\{x_{n}\} \in \overline{R(A)}$ with $\|x_{n}\|_{A} = 1$ such that $\|(\lambda -(T^{\diamond})^{\diamond})x_{n}\|_{A} \to 0$ as $n \to \infty$. Finally the equivalence extends to $\lambda \in \sigma_{A_{app}}((T^{\diamond})^{\diamond})$.
\end{proof}
In the following theorem, we prove the modified version of Theorem 2.11\cite{arxiv 007}.

\begin{theorem}\label{thm 8}
Let $T \in B_{A^{\frac{1}{2}}}(H)$. Then the boundary of $\sigma_{A}(T)$ is contained in the union of two sets $\sigma_{A_{app}}(T) \text{ and } (\sigma_{A_{app}}(T^{\diamond}))^{*} $, that is, $\partial {\sigma_{A}(T)} \subset \sigma_{A_{app}}(T) \cup (\sigma_{A_{app}}(T^{\diamond}))^{*}$. Here, $\partial {\sigma_{A}(T)}$ denotes the set of all boundary points.
\end{theorem}
\begin{proof}
Let us consider $\lambda \in \partial {\sigma_{A}(T)}$.  Then there exists a sequence $\{\lambda_{n}\}$ in $\mathbb {C} \setminus \sigma_{A}(T)$ such that $\lambda_{n} \to \lambda$ as $n \to \infty$. Moreover, we get $S_{\lambda_{n}} \in B_{A^{\frac{1}{2}}}(H) , \text { for all } n\in \mathbb N$, with $P (\lambda_{n} - T) S_{\lambda_{n}} = P S_{\lambda_{n}} (\lambda_{n} - T) = P$. Now we claim $\|S_{\lambda_{n}}\|_{A} \to \infty \text{ or } \|S^{\diamond}_{\lambda_{n}}\|_{A} \to \infty $ as $n \to \infty$. Suppose $\|S_{\lambda_{n}}\|_{A} \leq M$ and $\|S^{\diamond}_{\lambda_{n}}\|_{A} \leq M $,   for some $M > 0$ and for all $n \in \mathbb {N}$. We get a natural number $m\in \mathbb N$ such that $| \lambda_{m}- \lambda| < \frac{1}{M} \leq \frac{1}{{\|S_{\lambda_{m}}}\|_{A}}$ and  $| \lambda_{m}- \lambda| < \frac{1}{M} \leq \frac{1}{{\|S^{\diamond}_{\lambda_{m}}}\|_{A}}$ which imply $\|(\lambda_{m} - \lambda)S_{\lambda_{m}}\|_{A} < 1$ and $\|(\lambda_{m} - \lambda)S^{\diamond}_{\lambda_{m}}\|_{A} < 1$. From Lemma \ref{lemma 1}, we get that $(\lambda -T)$ is invertible in $B_{A^{\frac{1}{2}}}(H)$, which is a contradiction. Hence, $\|S_{\lambda_{n}}\|_{A} \to \infty$ or $\|S^{\diamond}_{\lambda_{n}}\|_{A} \to \infty$ as $n \to \infty$. Suppose $\|S_{\lambda_{n}}\|_{A} \to \infty$ as $n \to \infty$. So, we obtain a sequence $\{x_{n}\}$ in $H$, satisfying $\|x_{n}\|_{A} = 1$, such that $\alpha_{n} = \|S_{\lambda_{n}} x_{n}\|_{A} \geq \|S_{\lambda_{n}}\|_{A} - \frac{1}{n}, \text { for all } n\in \mathbb N$. Thus, $\alpha_{n} \to \infty$ as $n \to \infty$. Now, we consider elements $y_{n} = \frac{ S_{\lambda_{n}} x_{n}}{\alpha_{n}} \in H, \text { for all } n\in \mathbb N$. So, 
\begin{equation}\label{ineq 3}
\|(\lambda -T)y_{n}\|_{A} = \|(\lambda_{n} -T)y_{n} - (\lambda_{n} - \lambda)y_{n}\|_{A} \leq \frac{\|x_{n}\|_{A}}{\alpha_{n}} + |\lambda_{n}- \lambda|,      (\|y_{n}\|_{A} = 1,  \text { for all } n\in \mathbb N).
\end{equation}
The right-hand side of the above inequality (\ref{ineq 3}) goes to 0 as $n \to \infty$. This implies $\lambda \in \sigma_{A_{app}}(T)$. Similarly, we get $\lambda \in (\sigma_{A_{app}}(T^{\diamond}))^{*}$ when we consider $\|S^{\diamond}_{\lambda_{n}}\|_{A} \to \infty$ as $n \to \infty$. Therefore, $\partial {\sigma_{A}(T)} \subset \sigma_{A_{app}}(T) \cup (\sigma_{A_{app}}(T^{\diamond}))^{*}$.
\end{proof}

\begin{remark}\label{remark .24}
Let $T \in B_{A^{\frac{1}{2}}}(H)$ with $R(A)$ is closed or $ T \text{ commutes with } A$. Then $\sup\{\lambda : \lambda \in \sigma_{A}(H)\} = r_{A}(T)$ \cite{arxiv 001} which implies that Lemma \ref{lemma 1} is true without assuming $\|(T^{\prime}S)^{\diamond}\|_{A}< 1$. Moreover, we can state Theorem \ref{thm 8} as $\partial {\sigma_{A}(T)} \subset \sigma_{A_{app}}(T)$. The convex hull of $\sigma_{A}(T)$ is in $\overline {W_{A}(T)}$ because of Theorem \ref{thm 2.14} and Theorem \ref{thm 3.32}.
\end{remark}
\begin{theorem}\label{thm .24}
If $T$ belongs to the space $B_{A^{\frac{1}{2}}}(H)$, then the $A$-spectrum $\sigma_{A}(T)$ is non-empty.
\end{theorem}
\begin{proof}
The proof follows from Theorem 3.17 \cite{arxiv 001} and Corollary 5.7 \cite{MR4362420}.
\end{proof}

\begin{corollary}
Let $T \in B_{A^{\frac{1}{2}}}(H)$. Then either $\sigma_{A_{app}}(T)$ or $\sigma_{A_{app}}(T^{\diamond})$  is not empty.
\end{corollary}
\begin{proof}
We know that  $\sigma_{A}(T)$ is non-empty by Theorem \ref{thm .24}. The closedness of $\sigma_{A}(T)$ and Theorem \ref{thm 8} conclude that either $\sigma_{A_{app}}(T)$ or $\sigma_{A_{app}}(T^{\diamond})$ is a non-empty set in $\mathbb C$.
\end{proof}

\begin{theorem}\label{thm 30}
Let $T, T^{\prime} \in B_{A}(H)$ with $T^{\sharp} (T^{\prime})^{\sharp} = (T^{\prime})^{\sharp}T^{\sharp}$. Then the following conditions are equivalent:
\begin{enumerate}
    \item $T$ and $T^{\prime}$ both are $A$-invertible in $B_{A}(H)$.
    \item $TT^{\prime}$ is also $A$-invertible in $B_{A}(H)$.
\end{enumerate}
\end{theorem}
\begin{proof}
$(1)\implies (2)$ It can be directly shown by Proposition 4.5 in \cite{MR4362420} and Theorem \ref{thm 2}.

~ ~ ~~$(2) \implies (1)$ Let $S$ be an $A$-inverse in $B_{A}(H)$ of the $A$-invertible operator $(TT^{\prime})^{\sharp}$. Then $A(T^{\prime})^{\sharp} T^{\sharp}S = AS(T^{\prime})^{\sharp} T^{\sharp} = A$. We claim $A T^{\sharp}S (T^{\prime})^{\sharp} = A$. Let  $x \in N(A)$ be an arbitrary element.  Then $(T^{\prime})^{\sharp}x \in N(A)$ which implies $N(A) \subset N(A(T^{\prime})^{\sharp})$. Taking an element $z \in N(A(T^{\prime})^{\sharp})$ implies $P(T^{\prime})^{\sharp}z = 0$. So, $Pz= PS(T^{\prime})^{\sharp}T^{\sharp}z= PST^{\sharp}P(T^{\prime})^{\sharp}z = 0$. Hence, the reverse inclusion is also true and $N(A) = N(A(T^{\prime})^{\sharp})$. So, $ A(T^{\prime})^{\sharp} T^{\sharp}S =A$
implies $A(T^{\prime})^{\sharp} T^{\sharp}S(T^{\prime})^{\sharp} = A(T^{\prime})^{\sharp}$ and
     \begin{align*}\label{equ 11}
                 A(T^{\prime})^{\sharp}T^{\sharp}S(T^{\prime})^{\sharp}=A(T^{\prime})^{\sharp}\\
                A(T^{\prime})^{\sharp}(T^{\sharp}S(T^{\prime})^{\sharp} - I)= 0\\
                AT^{\sharp}S(T^{\prime})^{\sharp} = A.
     \end{align*}
Thus, $(T^{\prime})^{\sharp}$ is $A$-invertible with an $A$-inverse $T^{\sharp}S \in B_{A}(H)$. Similarly, we can prove that $T^{\sharp}$ is also $A$-invertible in $B_{A}(H)$. Therefore, $T \text{ and } T^{\prime}$ both are $A$-invertible in $B_{A}(H)$.
\end{proof} 
\begin{theorem}\label{thm 3.34}
Let $T \in B_{A}(H)$. Then $T$ is $A$-selfadjoint if and only if $W_{A}(T)$ is real.
\end{theorem}
\begin{proof}
Since $T^{*}A = AT$, we have
\begin{align*}
\langle Tx, x \rangle_{A} = \langle T^{*}Ax, x\rangle = \langle x, ATx\rangle = \overline{\langle Tx, x \rangle_{A}} , \text { for all } x \in H.
\end{align*}
Thus, $W_{A}(T)$ is real.

Conversely, for all $x\in H$ with $\|x\|_{A} = 1$, $
\langle Tx, x \rangle_{A} = \overline{\langle Tx, x \rangle_{A}}$. Hence $\langle (T^{*}A -AT)x, x\rangle = 0.$ Let us consider $z \in \overline{R(A)}$. Then $\langle (T^{*}A - AT)\frac{z}{\|z\|_{A}}, \frac{z}{\|z\|_{A}}\rangle = 0$ implies $\langle (T^{*}A - AT)z, z\rangle = 0$. So, $T^{*}A = AT$ in $\overline{R(A)}$. Moreover, $T^{*}A = AT$ in $N(A)$ because $T(N(A)) \subset N(A)$. Therefore, $T$ is $A$-selfadjoint. 
\end{proof}
\begin{corollary}
Let $T \in B_{A^{\frac{1}{2}}}(H)$ with $R(A)$ is closed. If $T$ is $A$-selfadjoint, then $\sigma_{A}(T)$ is real.
\end{corollary}
\begin{proof}
Drawing upon Remark \ref{remark .24} and Theorem \ref{thm 3.34}, it can be inferred that $\sigma_{A}(T)$ is real.
\end{proof}
\begin{theorem}
 If $T$ belongs to the space $B_{A}(H)$ and is  $A$-invertible in $B_{A}(H)$ with an $A$-inverse $S$, then $T$ is $A$-selfadjoint if and only if $S$ is also $A$-selfadjoint.
 \end{theorem}
 \begin{proof}
 By the given conditions, we have $T^{*}A = AT$ and $AST = ATS =A$. Then $S^{*}T^{*}A = A$ if and only if $S^{*} AT = A$. So, $S^{*}ATS = AS$ implies $S^{*}A = AS$. Therefore, $S$ is $A$-selfadjoint.
 To prove the converse part, we interchange $S$ and $T$, respectively. 
 \end{proof}

\section{$A$-approximate point spectrum of tensor product of two $A^{\frac{1}{2}}$-adjoint operators}

 For $x\in H_{1}$ and $y\in H_{2}$, we define a linear operator $x \otimes y$ from $H_{2}$ to $H_{1}$ by 
\begin{align*}
(x \otimes y)(z) = \langle z, \overline y\rangle x,  \text { for all } z\in H_{2}.
\end{align*}
The expressions of the form $x \otimes y$ are commonly referred to as ``elementary tensors". Let $H_{1} \odot H_{2}$ be the linear space spanned by elementary tensors $x \otimes y$. For $u = \sum_{k=1}^{r} c_{k}(x_{k} \otimes y_{k}) \in H_{1} \odot H_{2}$ and $v = \sum_{\ell=1}^{s} d_{\ell}({x_{\ell}}^{'} \otimes {y_{\ell}}^{'}) \in H_{1} \odot H_{2}$, we define $(u,v) = \sum_{k=1}^{r} \sum_{\ell=1}^{s} c_{k}{\overline d_{\ell}} \langle x_{k}, {x_{\ell}}^{'}\rangle \langle y_{k}, {y_{\ell}}^{'}\rangle$. It is easy to show that $(H_{1} \odot H_{2}, (.,.))$ is a well-defined inner product space. Now, the completion of the inner product space ($H_{1} \odot H_{2}, (.,.))$ is called the tensor product  of Hilbert spaces $H_{1}, H_{2}$, denoted by  $H_{1} \otimes H_{2}$. The norm in the Hilbert space $H_{1} \otimes H_{2}$ has the cross-property 
\begin{align*}
\|x \otimes y\| = \|x\| \|y\|,    \text { for all } x\in H_{1},  \text { for all } y\in H_{2}.
\end{align*}
Let $A_{1}$ and $A_{2}$ be two positive operators in two Hilbert spaces $H_{1}$ and $H_{2}$ respectively with $T_{1}\in B_{A_{1}^{\frac{1}{2}}}(H_{1}), T_{2} \in B_{A_{2}^{\frac{1}{2}}}(H_{2})$. We now define $T_{1} \odot T_{2}$ on $H_{1} \odot H_{2}$ by $T_{1} \odot T_{2} (a_{k}\sum_{k=1}^{m} z_{k} \otimes w_{k}) = \sum_{k=1}^{m} a_{k}(T_{1}z_{k} \otimes T_{2}w_{k})$. Then $T_{1} \odot T_{2}$ is bounded. We extend the bounded operator $T_{1} \odot T_{2}$ to Hilbert space $H_{1} \otimes H_{2}$ and the extended operator is called the tensor product of two operators $T_{1}$ and $T_{2}$, denoted by $T_{1} \otimes T_{2}$. It can be easily shown that $T_{1} \otimes T_{2}$ is in $B_{{A_{1}^{\frac{1}{2}}} \otimes A_{2}^{\frac{1}{2}}}(H_{1} \otimes H_{2})$.

In this section, we discuss the $A$-approximate point spectrum of the tensor product of two $A^{\frac{1}{2}}$-adjoint operators.

\begin{theorem} Let $A_{1}$ and $A_{2}$ be two positive operators in two Hilbert spaces $H_{1}$ and $H_{2}$ respectively with $T_{1}\in B_{A_{1}^{\frac{1}{2}}}(H_{1}), T_{2} \in B_{A_{2}^{\frac{1}{2}}}(H_{2})$. Then $$\sigma_{{A_{1}}_{app}}(T_{1})  \sigma_{{A_{2}}_{app}}(T_{2}) \subset \sigma_{({A_{1} \otimes A_{2})}_{app}}(T_{1} \otimes T_{2}).$$
\end{theorem}
\begin{proof}
Let us consider $\lambda_{1} \in \sigma_{{A_{1}}_{app}}(T_{1})$ and $\lambda_{2} \in \sigma_{{A_{2}}_{app}}(T_{2})$. Then there exist  sequences $\{x_{n}\}$ and $\{y_{n}\}$ in $H_{1}$ and $H_{2}$ respectively with $\|x_{n}\|_{A_{1}} = \|y_{n}\|_{A_{2}} = 1$ such that $\|(\lambda_{1} - T_{1})x_{n}\|_{A_{1}} \to 0$ and $\|(\lambda_{2} - T_{2})y_{n}\|_{A_{2}} \to 0$ as $n \to \infty$. Moreover, $\|x_{n} \otimes y_{n}\|_{A_{1} \otimes A_{2}} = \|x_{n}\|_{A_{1}} \|y_{n}\|_{A_{2}} =1,  \text { for all } n\in \mathbb N$.
\begin{equation}\label{ineq 14}
\begin{split}
&~~~~\|((T_{1} \otimes T_{2}) - \lambda_{1} \lambda_{2})(x_{n} \otimes y_{n})\|_{A_{1} \otimes A_{2}}\\
&= \|{(A_{1}}^{\frac{1}{2}}T_{1}x_{n} \otimes {A_{2}}^{\frac{1}{2}} T_{2}y_{n})- \lambda_{1} \lambda_{2} ({A_{1}}^{\frac{1}{2}} x_{n} \otimes {A_{2}}^{\frac{1}{2}}y_{n})\|\\
&=\|({A_{1}}^{\frac{1}{2}}(T_{1}- \lambda_{1}) \otimes {A_{2}}^{\frac{1}{2}}T_{2} + \lambda_{1}{A_{1}}^{\frac{1}{2}} \otimes {A_{2}}^{\frac{1}{2}}(T_{2} - \lambda_{2}))(x_{n} \otimes y_{n})\|\\
&\leq \|(T_{1}- \lambda_{1})x_{n}\|_{A_{1}} \|T_{2}y_{n}\|_{A_{2}} + |\lambda_{1}| \|(T_{2} -\lambda_{2})y_{n}\|_{A_{2}}.
\end{split}
\end{equation}
The right-hand side of the inequality (\ref{ineq 14}) goes to $0$ as $n \to \infty$. So, $\lambda_{1} \lambda_{2} \in \sigma_{({A_{1} \otimes A_{2})}_{app}}(T_{1} \otimes T_{2})$. Therefore, $\sigma_{{A_{1}}_{app}}(T_{1})  \sigma_{{A_{2}}_{app}}(T_{2}) \subset \sigma_{({A_{1} \otimes A_{2})}_{app}}(T_{1} \otimes T_{2})$.
\end{proof}
\begin{lemma} Let $A_{1}$ and $A_{2}$ be two positive operators in  Hilbert spaces $H_{1}$ and $H_{2}$ respectively with $T_{1}\in B_{A_{1}^{\frac{1}{2}}}(H_{1}), T_{2} \in B_{A_{2}^{\frac{1}{2}}}(H_{2})$. Then $$\sigma_{A_{1} \otimes A_{2}}(T_{1} \otimes I) \subset \sigma_{A_{1}}(T_{1}) \ \ \text{and} \ \ \sigma_{A_{1} \otimes A_{2}}(I \otimes T_{2}) \subset \sigma_{A_{2}}(T_{2}).$$
\end{lemma}
\begin{proof}
Let us consider $\lambda \in \rho_{A_{1}}(T_{1})$. Then there exists $S_{\lambda} \in B_{{A_{1}}^{\frac{1}{2}}}(H)$ such that $A_{1}(\lambda -T_{1}) S_{\lambda} = A_{1}S_{\lambda} (\lambda -T_{1}) = A_{1}$. So,
\begin{align*}
(A_{1} \otimes A_{2}) ((\lambda I \otimes I) - (T_{1} \otimes I)) (S_{\lambda} \otimes I) &= A_{1} \otimes A_{2}\\ &=(A_{1} \otimes A_{2})(S_{\lambda} \otimes I)  ((\lambda I \otimes I) - (T_{1} \otimes I)).
\end{align*}
Thus, $\lambda \in \rho_{A_{1} \otimes A_{2}}(T_{1} \otimes I)$. Therefore, $\sigma_{A_{1} \otimes A_{2}}(T_{1} \otimes I) \subset \sigma_{A_{1}}(T_{1})$. Similarly, we can prove that $\sigma_{A_{1} \otimes A_{2}}(I \otimes T_{2}) \subset \sigma_{A_{2}}(T_{2})$.
\end{proof}
\noindent We now consider a special case when  $A_{1}= A_{2} = A$ and $H_{1} = H_{2}= H$.
\begin{theorem} Let $A$  be a positive operator in a Hilbert space $H$ with $T_{1}\in B_{A^{\frac{1}{2}}}(H), T_{2} \in B_{A^{\frac{1}{2}}}(H)$. Then 
$$\sigma_{{(A \otimes A)}_{app}} (T_{1} \odot I) = \sigma_{A_{app}}(T_{1}) \ \ \text{and} \ \ \sigma_{{(A \otimes A)}_{app}} (I \odot T_{2}) = \sigma_{A_{app}}(T_{2}).$$
\end{theorem}
\begin{proof}
Let $\lambda \in \sigma_{{(A \otimes A)}_{app}} (T_{1} \odot I)$. Then we obtain a sequence  $\{x_{n} \otimes y_{n}\}$ in $H \odot H$, where $\|x_{n} \otimes y_{n}\|_{A \otimes A} = 1$, such that $\|(\lambda(I \odot I) - (T_{1} \odot I)) x_{n} \otimes y_{n}\|_{A \otimes A} \to 0$ as $n \to \infty$. Moreover, $\|x_{n} \otimes y_{n}\|_{A \otimes A} = \|x_{n}\|_{A} \|y_{n}\|_{A} =1, \text { for all } n\in \mathbb N$. Assume, $\lambda \notin \sigma_{A_{app}}(T_{1})$. Then there exists $d_{\lambda} > 0$ such that $\|(\lambda -T_{1})x\|_{A} \geq d_{\lambda} \|x\|_{A}, \text{ for all } x\in H$. So,
\begin{align}\label{ineq 15}
d_{\lambda} = d_{\lambda} \|x_{n}\|_{A} \|y_{n}\|_{A} \leq \|(\lambda -T_{1})x_{n}\|_{A} \|y_{n}\|_{A} = \|(\lambda -T_{1})x_{n} \otimes y_{n}\|_{A \otimes A}.
\end{align}
The right-hand side of the inequality (\ref{ineq 15}) goes to $0$ as $n \to \infty$. So $d_{\lambda} = 0$, which is contradiction. Hence,  $\sigma_{{(A \otimes A)}_{app}} (T_{1} \odot I) \subset \sigma_{A_{app}}(T_{1})$.

Conversely, let $\mu \in \sigma_{A_{app}}(T_{1})$. We obtain a sequence $\{z_{n}\}$ in $H$ with $\|z_{n}\|_{A} =1$ such that $\|(\mu -T_{1}) z_{n}\|_{A} \to 0$ as $n \to \infty$. Moreover, $\|z_{n} \otimes z_{n}\|_{A \otimes A} = \|z_{n}\|_{A} \|z_{n}\|_{A} = 1$.
\begin{align}\label{eq 16}
\| (\mu (I \odot I) - (T_{1} \odot I)) (z_{n} \otimes z_{n})\|_{A \otimes A} = \|(\mu - T_{1})z_{n}\|_{A}.
\end{align}
The right-hand side of the equation (\ref{eq 16}) goes to $0$ as $n \to \infty$. Therefore, $\sigma_{A_{app}}(T_{1}) = \sigma_{{(A \otimes A)}_{app}} (T_{1} \odot I)$.
Similarly, we can prove that $\sigma_{{(A \otimes A)}_{app}} (I \odot T_{2}) = \sigma_{A_{app}}(T_{2})$.
\end{proof}
\section{Some characterizations of $A$-normal operators}
In the paper \cite{MR2878497}, an $A$-normal operator $T \in B_{A}(H)$ is defined as an operator for which $T^{\sharp} T = T T^{\sharp}$. This condition is equivalent to asserting that the range of $TT^{\sharp}$ is contained in the closure of the range of the operator $A$, and  $\|Tx\|_{A} = \|T^{\sharp}x\|_{A}, \text{ for all } x\in H$. 

On the other hand, the paper cited as \cite{arxiv 007} introduces a different definition of an $A$-normal operator $T$ in the space $B_{A}(H)$. In this case an operator is considered $A$-normal if it satisfies the condition $\|Tx\|_{A} = \|T^{\sharp}x\|_{A}, \text{ for all } x\in H$. This definition is used throughout in the paper. 

Additionally the concept of an $A$-hyponormal operator is introduced in the context of the paper \cite{MR3513410}. An operator $T \in B_{A}(H)$ is said to be $A$-hyponormal if $T^{\sharp}T - TT^{\sharp} \geq_{A} 0$, which is further equivalent to the condition $\|Tx\|_{A} \geq \|T^{\sharp}x\|_{A}, \text{ for all } x\in H$.

We introduce a new operator $T_{a}: \overline{R(A)} \mapsto \overline{R(A)}$ defined as follows: $\text{ for all } x\in \overline{R(A)}$ $T_{a}(x)= (A^{\frac{1}{2}})^{\dagger} T^{*}A^{\frac{1}{2}}(x)$, where $T \in B_{A^{\frac{1}{2}}}(H)$.
\begin{theorem}\label{thm 0.5.1}
Let $T \in B_{A}(H)$. Then $T$ is a $A$-normal operator if and only if $T_{a}$ is normal.
\end{theorem}
\begin{proof}
Suppose that $T$ is $A$-normal. It is easy to show that $T_{a}^{*} = {\overline{A^{\frac{1}{2}}T (A^{\frac{1}{2}})^{\dagger}}\vert}_{\overline{R(A)}}$. Now consider an arbitrary element $Aw \in R(A)$, we get
\begin{align*}
\|T_{a}(Aw)\| = \|T^{\sharp} A^{\frac{1}{2}}w\|_{A}= \|T A^{\frac{1}{2}}w\|_{A}= \|A^{\frac{1}{2}} T (A^{\frac{1}{2}})^{\dagger} Aw\|= \|T_{a}^{*} Aw\|.
\end{align*}
The denseness of $R(A)$ in $\overline{R(A)}$ confirms  $T_{a}$ is normal.

Conversely, let us consider $x \in H$. Then 
\begin{align*}
&\|T_{a}(Ax)\|= \|T_{a}^{*}(Ax)\| \\
& \|T^{\sharp}A^{\frac{1}{2}}x\|_{A} = \|T A^{\frac{1}{2}}x\|_{A}.\\
\end{align*}
We know that $R(A^{\frac{1}{2}})$ is dense in $\overline{R(A)}$. So, $\|Tu\|_{A}= \|T^{\sharp}u\|_{A}, \text{ for all } u \in \overline{R(A)}$. Again, $\|Tv\|_{A} = \|T^{\sharp}v\|_{A}, \text{ for all }v\in N(A)$. Therefore, $T$ is $A$-normal.
\end{proof}

\begin{theorem}
Let $T \in B_{A}(H)$. Then $T$ is $A$-selfadjoint if and only if $\tilde{T}$ is also selfadjoint.
\end{theorem}
\begin{proof}
First, we consider $T^{*}A = AT$. We will show that $\tilde {T}$ is self-adjoint. 
\begin{align*}
&(\tilde{T}Ax, Ay)\\ 
&= (\tilde{T}W_{A}x, Ay)\\
&= \langle A^{\frac{1}{2}}Tx, A^{\frac{1}{2}}y \rangle \\
&= \langle x, T^{*}Ay \rangle \\
&= \langle A^{\frac{1}{2}}x, A^{\frac{1}{2}}Ty \rangle \\
&= (Ax, W_{A}Ty)\\
&= (Ax, \tilde{T}Ay).
\end{align*}
The range of $A$ is dense in $\mathbf{R}(A^{\frac{1}{2}})$. Thus, $(\tilde{T})^{*} = \tilde{T}$.

Conversely, assume that $\tilde{T}$ is selfadjoint. Then for all $x, y \in H$, we get
\begin{align*}
& (\tilde{T}Ax, Ay) = (Ax, \tilde{T}Ay)\\
& (W_{A}Tx, Ay) = (Ax, W_{A}Ty)\\
& \langle A^{\frac{1}{2}}Tx, A^{\frac{1}{2}}y \rangle = \langle A^{\frac{1}{2}}x, A^{\frac{1}{2}}Ty \rangle \\
& \langle x, T^{*}Ay \rangle = \langle x, ATy \rangle.
\end{align*}
Therefore, $T$ is a $A$-selfadjoint operator.
\end{proof}
\begin{theorem}\cite{arxiv 007}
Let $T \in B_{A}(H)$. If $T$ is $A$-normal, then the following statements hold:
\begin{enumerate}
\item $(T- \lambda)$ is $A$-normal, for all $\lambda \in \mathbb C$.
\item The spectral radius of $T$ is equal to $\|T\|_{A}$.
\end{enumerate}
\end{theorem}
\begin{theorem}\label{thm 0.5.4}
Let $T \in B_{A}(H)$ be $A$-normal with $R(A)$ is closed. Then $\sigma_{A}(T) = \sigma_{A_{app}}(T)$.
\end{theorem}
\begin{proof}
Let us assume $\lambda \in \sigma_{A}(T) \setminus \sigma_{A_{app}}(T)$. Then there exists $d_{\lambda} > 0$ such that $\|(T - \lambda)x\|_{A} \geq d_{\lambda} \|x\|_{A}, \text{ for all }x\in H$. So, $A$-normality of $(T-\lambda)$ confirms that  $\|(T - \lambda)^{\sharp}x\|_{A} \geq d_{\lambda} \|x\|_{A}, \text{ for all }x\in H$. Again, $(A^{\frac{1}{2}})^{\dagger}$ is bounded because $R(A)$ is closed.
\begin{align*}
& \|(T - \lambda)^{\sharp}x\|_{A}^{2} \geq d_{\lambda}^{2} \|x\|_{A}^{2}\\
& \langle A(T- \lambda)A^{\dagger}(T-\lambda)^{*}Ax,x\rangle \geq d_{\lambda}^{2} \langle Ax,x \rangle.
\end{align*}
Then, there exists a reduced solution $V \in B(H)$ such that
\begin{align*}
& A(T- \lambda)(A^{\frac{1}{2}})^{\dagger} V = A^{\frac{1}{2}}\implies A(T- \lambda)(A^{\frac{1}{2}})^{\dagger} V A^{\frac{1}{2}}= A.\\
\end{align*}
Moreover, $(A^{\frac{1}{2}})^{\dagger} V A^{\frac{1}{2}} \in B_{A^{\frac{1}{2}}}(H)$.  The condition  $\|(T - \lambda)x\|_{A} \geq d_{\lambda} \|x\|_{A}, \text{ for all }x\in H$, confirms that $ N(A(T-\lambda)) \subset N(A)$. Now, the relation $A(T- \lambda) = ((T-\lambda)^{\sharp})^{*} A$ guarantees that $ N(A(T-\lambda))$ is equal to $N(A)$. By Proposition 3.7 \cite{MR4362420}, we get $\lambda \in \rho_{A}(T)$ which is a contradiction. So, $\lambda \in \sigma_{A_{app}}(T)$. Therefore, $\sigma_{A}(T)=\sigma_{A_{app}}(T)$.
\end{proof}
From the Spectral Mapping Theorem, we know that $\sigma(f(T)) = f(\sigma(T)), \text{ where } T$ is a normal operator and $f$ is a continuous function on $\sigma(T)$. It is obvious to ask whether the set $g(\sigma_{A}(T))$ is equal or not with the $A$-spectrum of a $B_{A^{\frac{1}{2}}}(H)$ operator, $\text{ where } g$ is a continuous function on $\sigma_{A}(T)$ and $T$ is $A$-normal. In the following theorems, our main goal is to show $g(\sigma_{A}(T)) = \sigma_{A}((A^{\frac{1}{2}})^{\dagger} g(T_{a}^{*})A^{\frac{1}{2}})$, when $R(A)$ is assumed to be closed.
\begin{proposition}\label{pro 0.5.5}
Let $T\in B_{A}(H)$ and $q(x,y)= a_{0} + a_{1}x^{r_{1}}y^{s_{1}} + \dots + a_{k}x^{r_{k}}y^{s_{k}}, \text{ where } a_{i} \in \mathbb{C}, i= 1, 2, \dots, k$. Then $\|q(T, T^{\sharp})\|_{A} = \|q(T_{a}^{*}, T_{a})\|$.
\end{proposition}
\begin{proof}
It is evident that $q(T, T^{\sharp})$  belongs to the space $B_{A}(H)$. Now,
\begin{align*}
&\|q(T,T^{\sharp})\|_{A} \\
&= \|\overline{A^{\frac{1}{2}}(a_{0} + a_{1}T^{r_{1}}(T^{\sharp})^{s_{1}} + \dots + a_{k}T^{r_{k}}(T^{\sharp})^{s_{k}}) (A^{\frac{1}{2}})^{\dagger}}\|\\
&= \|a_{0}\overline{A^{\frac{1}{2}} (A^{\frac{1}{2}})^{\dagger}} + a_{1} \overline{A^{\frac{1}{2}} T^{r_{1}}(T^{\sharp})^{s_{1}} (A^{\frac{1}{2}})^{\dagger}} +\dots + a_{k} \overline{A^{\frac{1}{2}} T^{r_{k}}(T^{\sharp})^{s_{k}} (A^{\frac{1}{2}})^{\dagger}}\|\\
&= \|a_{0}\overline{A^{\frac{1}{2}} (A^{\frac{1}{2}})^{\dagger}} + a_{1}(\overline{A^{\frac{1}{2}}T(A^{\frac{1}{2}})^{\dagger}})^{r_{1}} (\overline{A^{\frac{1}{2}}T^{\sharp}(A^{\frac{1}{2}})^{\dagger}})^{s_{1}} + \dots +  a_{k}(\overline{A^{\frac{1}{2}}T(A^{\frac{1}{2}})^{\dagger}})^{r_{k}} (\overline{A^{\frac{1}{2}}T^{\sharp}(A^{\frac{1}{2}})^{\dagger}})^{s_{k}}\|\\
&= \|a_{0}+ a_{1}(T_{a}^{*})^{r_{1}}(T_{a})^{s_{1}} +\dots + a_{k}(T_{a}^{*})^{r_{k}}(T_{a})^{s_{k}}\|\\
&= \|q(T_{a}^{*}, T_{a})\|.
\end{align*}
\end{proof}
\begin{theorem}\label{thm 0.5.6}
Let $T \in B_{A}(H)$ be $A$-normal with $R(A)$ is closed. Then $\sigma_{A}(T) = \sigma(T_{a}^{*})$. 
\end{theorem}
\begin{proof}
Theorem \ref{thm 0.5.1} says that $T_{a}^{*}$ is normal. So, $\sigma(T_{a}^{*}) = \sigma_{app}(T_{a}^{*})$, where $\sigma_{app}(T_{a}^{*})$ is the approximate point spectrum of $T_{a}^{*}$. Now, take any arbitrary $\lambda \in \rho_{A}(T)$, and there exists an operator $S_{\lambda} \in B_{A}(H)$ such that 
\begin{equation}\label{eqn 0.14}
P(\lambda -T) S_{\lambda} = P S_{\lambda}(\lambda -T) = P.
\end{equation}
From the equality (\ref{eqn 0.14}), we get
\begin{align*}
&S_{\lambda}^{*}(\lambda - T)^{*}P = (\lambda - T)^{*}S_{\lambda}^{*}P= P\\
&\implies (A^{\frac{1}{2}})^{\dagger}S_{\lambda}^{*} A^{\frac{1}{2}} (A^{\frac{1}{2}})^{\dagger} (\lambda-T)^{*} A^{\frac{1}{2}}= (A^{\frac{1}{2}})^{\dagger}(\lambda-T)^{*}A^{\frac{1}{2}} (A^{\frac{1}{2}})^{\dagger}S_{\lambda}^{*}A^{\frac{1}{2}}= P\\
&\implies (S_{\lambda})_{a}(\overline{\lambda} - T_{a}) = (\overline{\lambda} - T_{a})(S_{\lambda})_{a}= I_{R(A)}, \text{ where } (S_{\lambda})_{a}= (A^{\frac{1}{2}})^{\dagger}S_{\lambda}^{*} A^{\frac{1}{2}}\vert_{R(A)}.
\end{align*}
Thus, $\sigma(T_{a}^{*}) \subset \sigma_{A}(T)$. Now, we will show the reverse inclusion $\sigma_{A}(T) \subset \sigma(T_{a}^{*})$.

Let $\mu \in \sigma_{A}(T)$. From Theorem \ref{thm 0.5.4}, we get $\mu \in \sigma_{A_{app}}(T)$. Then, we obtain a sequence $\{z_{n}\}$ in $R(A)$ with $\|z_{n}\|_{A}=1$ such that $\|(\mu - T)z_{n}\|_{A} \to 0$, as $n \to \infty$. So, $\|(\mu - T_{a}^{*})A^{\frac{1}{2}}z_{n}\| = \|(\mu -T)z_{n}\|_{A} \to 0$ as $n \to \infty$. Hence, $\mu \in \sigma_{app}(T_{a}^{*})= \sigma(T_{a}^{*})$.
 Therefore, $\sigma_{A}(T) = \sigma(T_{a}^{*})$.
\end{proof}
\begin{theorem}
Let $T \in B_{A}(H)$ be $A$-normal with $R(A)$ be closed. Then $f(\sigma_{A}(T)) = \sigma_{A}((A^{\frac{1}{2}})^{\dagger} f(T_{a}^{*}) A^{\frac{1}{2}})$, where $f$ is a continuous function on $\sigma_{A}(T)$.
\end{theorem}
\begin{proof}
Let us consider $\lambda \in \rho(f(T_{a}^{*}))$. Then there exists an operator $V_{\lambda} \in B(R(A))$ such that
\begin{equation}\label{eqn 0.15}
(\lambda -f(T_{a}^{*}))V_{\lambda} = V_{\lambda}(\lambda - f(T_{a}^{*})) = I_{R(A)}.
\end{equation}
From the equality (\ref{eqn 0.15}), we get
\begin{align*}
&(A^{\frac{1}{2}})^{\dagger}(\lambda - f(T_{a}^{*}))A^{\frac{1}{2}} (A^{\frac{1}{2}})^{\dagger}V_{\lambda}A^{\frac{1}{2}}= (A^{\frac{1}{2}})^{\dagger}V_{\lambda} A^{\frac{1}{2}}(A^{\frac{1}{2}})^{\dagger}(\lambda - f(T_{a}^{*})) A^{\frac{1}{2}} = P\\
& P(\lambda - (A^{\frac{1}{2}})^{\dagger} f(T_{a}^{*})A^{\frac{1}{2}}) (A^{\frac{1}{2}})^{\dagger}V_{\lambda}A^{\frac{1}{2}} = P (A^{\frac{1}{2}})^{\dagger}V_{\lambda}A^{\frac{1}{2}}(\lambda - (A^{\frac{1}{2}})^{\dagger} f(T_{a}^{*})A^{\frac{1}{2}}) = P.
\end{align*}
It is easy to prove that $(A^{\frac{1}{2}})^{\dagger} V_{\lambda} A^{\frac{1}{2}} \in B_{A^{\frac{1}{2}}}(H)$. So, $\lambda \in \rho_{A}((A^{\frac{1}{2}})^{\dagger} f(T_{a}^{*}) A^{\frac{1}{2}})$. Thus, $\sigma_{A}((A^{\frac{1}{2}})^{\dagger} f(T_{a}^{*}) A^{\frac{1}{2}}) \subset \sigma(f(T_{a}^{*}))$. Now, we claim the reverse inclusion $\sigma(f(T_{a}^{*})) \subset \sigma_{A}((A^{\frac{1}{2}})^{\dagger} f(T_{a}^{*}) A^{\frac{1}{2}})$.

Let, $\mu \in \rho_{A}((A^{\frac{1}{2}})^{\dagger} f(T_{a}^{*}) A^{\frac{1}{2}})$. Then there exists $S_{\mu} \in B_{A^{\frac{1}{2}}}(H)$ such that
\begin{equation}\label{eqn 0.16}
P(\mu - ((A^{\frac{1}{2}})^{\dagger} f(T_{a}^{*}) A^{\frac{1}{2}}))S_{\mu} = P S_{\mu} ((\mu - ((A^{\frac{1}{2}})^{\dagger} f(T_{a}^{*}) A^{\frac{1}{2}})))=P.
\end{equation}
From the equality (\ref{eqn 0.16}), we get
\begin{align*}
& S_{\mu}^{*}(A^{\frac{1}{2}})(\overline{\mu} - (f(T_{a}^{*}))^{*}) (A^{\frac{1}{2}})^{\dagger} = A^{\frac{1}{2}} (\overline{\mu} - (f(T_{a}^{*}))^{*}) (A^{\frac{1}{2}})^{\dagger}S_{\mu}^{*}= I_{R(A)}\\
& \implies (A^{\frac{1}{2}})^{\dagger}S_{\mu}^{*}(A^{\frac{1}{2}}) (\overline{\mu} - (f(T_{a}^{*}))^{*})= (\overline{\mu} - (f(T_{a}^{*}))^{*}) (A^{\frac{1}{2}})^{\dagger}S_{\mu}^{*}(A^{\frac{1}{2}})= I_{R(A)}.
\end{align*}
So, $\mu \in \rho(f(T_{a}^{*}))$ which implies that $\sigma(f(T_{a}^{*})) \subset \sigma_{A}((A^{\frac{1}{2}})^{\dagger} f(T_{a}^{*}) A^{\frac{1}{2}})$. Therefore, $\sigma_{A}((A^{\frac{1}{2}})^{\dagger} f(T_{a}^{*}) A^{\frac{1}{2}})= \sigma(f(T_{a}^{*})) = f(\sigma(T_{a}^{*})) = f(\sigma_{A}(T))$.

\end{proof}

The closure of the numerical range of a normal operator coincides with the convex hull of its spectrum. A natural question arises: Is it true that $\overline{W_{A}(T)}= conv(\sigma_{A}(T))$? In \cite{arxiv 007}, authors claim the validity of this equality in Theorem 2.21 \cite{arxiv 007}. However, the proof of Theorem 2.21\cite{arxiv 007} is not rigorous, as the authors rely on Theorem 2.12 \cite{arxiv 007} to establish $conv(\sigma_{A}(T)) \subset \overline{W_{A}(T)}$. The proof of Theorem 2.12 \cite{arxiv 007} itself relies on Theorem 2.11 \cite{arxiv 007}, which is not precise due to the issues identified in Theorem 5.1\cite{MR4362420}. In this context, we will demonstrate the correctness of the equality  $\overline{W_{A}(T)}= conv(\sigma_{A}(T))$ when $T$ commutes with $A$ and $T$ is $A$-normal.

We introduce a new operator $T_{c}: \overline{R(A)} \mapsto \overline{R(A)}$ by $T_{c}(x) = \overline{A^{\frac{1}{2}}T (A^{\frac{1}{2}})^{\dagger}}\vert_{\overline {R(A)}}(x), \text{ for all } x\in \overline{R(A)}$ and $T \in B(A^{\frac{1}{2}}(H))$. 
\begin{lemma}\label{lemma 5.5}
Let $T \in B_{A^{\frac{1}{2}}}(H)$. Then $W_{A}(T) \subset W(T_{c})$.
\end{lemma}
\begin{proof}
Let $y \in W_{A}(T)$. Then we obtain $x \in \overline{R(A)}$ with $\|x\|_{A}= 1$ such that $y = \langle Tx, x \rangle_{A} = \langle T_{c}(A^{\frac{1}{2}}x), A^{\frac{1}{2}}x \rangle$. So, $y \in W(T_{c})$. Therefore, $W_{A}(T) \subset W(T_{c})$.
\end{proof}
\begin{theorem} \label{thm 5.6}
Let $T \in B_{A^{\frac{1}{2}}}(H)$. Then $\sigma(T_{c}) \subset \sigma_{A}(T)$.
\end{theorem}
\begin{proof}
Let $\lambda \in \rho_{A}(T)$. Then there exists $S_{\lambda} \in B_{A^{\frac{1}{2}}}(H)$ such that 
\begin{equation}\label{eqn 14}
A^{\frac{1}{2}}(\lambda -T) S_{\lambda} = A^{\frac{1}{2}} S_{\lambda}(\lambda -T) = A^{\frac{1}{2}}.
\end{equation}
Now we claim that $(\lambda -T_{c})$ is invertible with the inverse $\overline{A^{\frac{1}{2}} S_{\lambda}(A^{\frac{1}{2}})^{\dagger}}\vert_{\overline{R(A)}}$ . 

From equation (\ref{eqn 14}), we get 
\begin{align*}
& S_{\lambda}^{*}(\lambda-T)^{*}A^{\frac{1}{2}} = (\lambda -T)^{*}S_{\lambda}^{*}A^{\frac{1}{2}}= A^{\frac{1}{2}}\\
& \implies ((A^{\frac{1}{2}})^{\dagger} S_{\lambda}^{*} A^{\frac{1}{2}}) ((A^{\frac{1}{2}})^{\dagger} (\lambda -T)^{*}A^{\frac{1}{2}}) = ((A^{\frac{1}{2}})^{\dagger} (\lambda -T)^{*}A^{\frac{1}{2}})((A^{\frac{1}{2}})^{\dagger} S_{\lambda}^{*} A^{\frac{1}{2}})= P\\
&= ((A^{\frac{1}{2}})^{\dagger} (\lambda -T)^{*}A^{\frac{1}{2}})^{*} ((A^{\frac{1}{2}})^{\dagger} S_{\lambda}^{*} A^{\frac{1}{2}})^{*}= ((A^{\frac{1}{2}})^{\dagger} S_{\lambda}^{*} A^{\frac{1}{2}})^{*} ((A^{\frac{1}{2}})^{\dagger} (\lambda -T)^{*}A^{\frac{1}{2}})^{*} =P\\
&\implies (\overline{A^{\frac{1}{2}} (\lambda -T)(A^{\frac{1}{2}})^{\dagger}}\vert_{\overline{R(A)}}) (\overline{A^{\frac{1}{2}} S_{\lambda}(A^{\frac{1}{2}})^{\dagger}}\vert_{\overline{R(A)}}) = (\overline{A^{\frac{1}{2}} S_{\lambda}(A^{\frac{1}{2}})^{\dagger}}\vert_{\overline{R(A)}})  (\overline{A^{\frac{1}{2}} (\lambda -T)(A^{\frac{1}{2}})^{\dagger}}\vert_{\overline{R(A)}}) = I_{\overline{R(A)}}\\
&=(\lambda -T_{c}) (\overline{A^{\frac{1}{2}} S_{\lambda}(A^{\frac{1}{2}})^{\dagger}}\vert_{\overline{R(A)}}) = (\overline{A^{\frac{1}{2}} S_{\lambda}(A^{\frac{1}{2}})^{\dagger}}\vert_{\overline{R(A)}}) (\lambda - T_{c}) = I_{\overline{R(A)}}.
\end{align*}
Thus, $\lambda \in \rho(T_{c})$. Therefore, $\sigma(T_{c}) \subset \sigma_{A}(T)$.
\end{proof}
\begin{lemma}\label{lemma 5.7}
Let $T$ be $A$-normal. Then $T_{c}$ is normal.
\end{lemma}
\begin{proof}
Since $T$ is $A$-normal, $\|Tx\|_{A} = \|T^{\sharp}\|_{A}, \text{ for all } x\in H$. It is easy to prove that $T_{c}^{*} = (A^{\frac{1}{2}})^{\dagger} T^{*}A^{\frac{1}{2}}\vert_{\overline{R(A)}}$. Now, $\|T_{c}A^{\frac{1}{2}}x\| = \|Tx\|_{A} = \|T^{\sharp}x\|_{A} = \|T_{c}^{*}A^{\frac{1}{2}}x\|$. We know that $R(A^{\frac{1}{2}})$ is dense in $\overline{R(A)}$. Therefore, $T_{c}$ is normal.
\end{proof}
\begin{theorem}\label{thm 5.8}
Let $T$ be a $A$-normal operator. Then $\overline{W_{A}(T)} \subset conv(\sigma_{A}(T))$.
\end{theorem}
\begin{proof}
 Lemma \ref{lemma 5.5}, Theorem \ref{thm 5.6}, and Lemma \ref{lemma 5.7} confirm that
 \begin{align*}
 \overline{W_{A}(T)} \subset \overline{W(T_{c})} = conv(\sigma(T_{c})) \subset conv(\sigma_{A}(T)).
 \end{align*}
\end{proof}
\begin{corollary}
Let $T$ be a $A$-normal operator which commutes with $A$. Then $\overline{W_{A}(T)} = conv(\sigma_{A}(T))$.
\end{corollary}
\begin{proof}
Theorem \ref{thm 5.8} and Remark \ref{remark .24} show that $\overline{W_{A}(T)} = conv(\sigma_{A}(T))$.
\end{proof}
\begin{theorem}
Let $T$ be a $A$-normal operator. Then $T^{\diamond}$ is a normal operator.
\end{theorem}
\begin{proof}
We know that $T^{\diamond} = (A^{\frac{1}{2}})^{\dagger} T^{*}A^{\frac{1}{2}}$ and $(T^{\diamond})^{*} = \overline{A^{\frac{1}{2}} T (A^{\frac{1}{2}})^{\dagger}}$.
Then, $\|(T^{\diamond})^{*}A^{\frac{1}{2}}x\| = \|Tx\|_{A}= \|T^{\sharp}x\|_{A} = \|T^{\diamond}A^{\frac{1}{2}}x\|$. $R(A^{\frac{1}{2}})$ is dense in $\overline{R(A)}$. So, $\|T^{\diamond}u\| = \|(T^{\diamond})^{*}u\|, \text{ for all }u \in \overline{R(A)}$. Moreover, $\text{ for all }v \in N(A)$ we have $\|T^{\diamond}v\| = \|(T^{\diamond})^{*}v\|$. Therefore, $T^{\diamond}$ is normal.
\end{proof}
\begin{theorem}
Let $T \in B_{A}(H)$. Then $T$ is $A$-hyponormal if and only if $\tilde{T}$ is hyponormal.
\end{theorem}
\begin{proof}
First, we will show that $(\tilde T)^{*}Ay = T^{*}Ay, \text{ for all } y \in \overline{R(A)}$. For all $Ax \in R(A)$, we get
\begin{align*}
(\tilde{T}(Ax), Ay) = (ATx, Ay)= \langle ATx, Py \rangle = \langle Ax, T^{\sharp}Py \rangle= (Ax, AT^{\sharp}y).
\end{align*}
We know that $R(A)$ is dense is $\mathbf{R}(A^{\frac{1}{2}})$. So, $(\tilde T)^{*}Ay = AT^{\sharp}y =T^{*}Ay, \text{ for all } y \in \overline{R(A)}$.
Now, $\|\tilde{T}Au\|_{\mathbf{R}(A^{\frac{1}{2}})} = \|ATu\|_{\mathbf{R}(A^{\frac{1}{2}})} = \|Tu\|_{A}$ and $\|(\tilde{T})^{*}Au\|_{\mathbf{R}(A^{\frac{1}{2}})} = \|T^{*}Au\|_{\mathbf{R}(A^{\frac{1}{2}})} = \|T^{\sharp}u\|_{A}$.

When $T$ is $A$-hyponormal, we have $\|(\tilde{T})^{*}Au\|_{\mathbf{R}(A^{\frac{1}{2}})}=\|T^{\sharp}u\|_{A} \leq \|Tu\|_{A} = \|\tilde{T}Au\|_{\mathbf{R}(A^{\frac{1}{2}})}$. The density property of $R(A)$ in ${\mathbf{R}(A^{\frac{1}{2}})}$ says that $\tilde{T}$ is hyponormal.

Conversely, $\tilde{T}$ is hyponormal means $\|T^{\sharp}v\|_{A} = \|(\tilde{T})^{*}Av\|_{\mathbf{R}(A^{\frac{1}{2}})} \leq \|\tilde{T}Av\|_{\mathbf{R}(A^{\frac{1}{2}})} = \|Tv\|_{A}, \text{ for all }v \in H$ which equals to $T$ is $A$-hyponormal.

\end{proof}
\begin{theorem}\cite{MR0187093}\label{thm 5.12}
Let $T$ be a hyponormal operator. Then $\overline{W_{A}(T)} = conv(\sigma(T))$.
\end{theorem}
\begin{theorem}
Let $T\in B_{A}(H)$ be $A$-hyponormal which commutes with $A$. Then $\overline{W_{A}(T)} = conv(\sigma_{A}(T))$.
\end{theorem}
\begin{proof}
It is easy to show that $T_{c}$ is also hyponormal when $T$ is $A$-hyponormal. By Lemma \ref{lemma 5.5}, Theorem \ref{thm 5.6}, and Theorem \ref{thm 5.12} give that $\overline{W_{A}(T)} \subset \overline{W(T_{c})} = conv(\sigma(T_{c})) \subset conv(\sigma_{A}(T))$.
Hence by Remark \ref{remark .24}, we get $\overline{W_{A}(T)} = conv(\sigma_{A}(T))$.
\end{proof}

\begin{center}
	\textbf{Acknowledgements}
\end{center}

\noindent The present work of the second author was partially supported by Science and Engineering Research Board (SERB), Department of Science and Technology, Government of India (Reference Number: MTR/2023/000471) under the scheme ``Mathematical Research Impact Centric Support (MATRICS)''.

\end{document}